\documentclass[12pt,a4paper]{amsart}
\usepackage{amssymb,amsxtra}
\usepackage[english,french]{babel}

\newtheorem{thm}{Theorem}[section]
\newtheorem{lem}[thm]{Lemma}

\theoremstyle{definition}

\theoremstyle{remark}

\theoremstyle{plain}

\theoremstyle{remark}

\newtheorem*{example}{Example}
\numberwithin{equation}{section}

\begin{document}

\title{ Enumeration of some particular   quintuple  Persymmetric  Matrices over $\mathbb{F}_{2} $ by rank}
\author{Jorgen~Cherly}
\address{D\'epartement de Math\'ematiques, Universit\'e de
    Brest, 29238 Brest cedex~3, France}
\email{Jorgen.Cherly@univ-brest.fr}
\email{andersen69@wanadoo.fr}

\maketitle 
\begin{abstract}
Dans cet article nous comptons le nombre de certaines  quintuples  matrices persym\' etriques de rang i sur $ \mathbb {F} _ {2} . $

 \end{abstract}

\selectlanguage{english}

\begin{abstract}
In this paper we count the number of some particular  quintuple  persymmetric rank i matrices over  $ \mathbb{F}_{2}.$
 \end{abstract}
 
  \maketitle 
\newpage
\tableofcontents
\newpage

  \section{Introduction}
  \label{sec 1}  
  In this paper we propose to compute in the most simple case  the number of  quintuple persymmetric 
  matrices with entries in   $ \mathbb{F}_{2}$ of rank i\\
  That is to compute the number  $ \Gamma_{i}^{\left[2\atop{2 \atop{ 2\atop {2\atop2}}}\right]\times k}$ 
  of  quintuple  persymmetric matrices in   $ \mathbb{F}_{2}$ of rank i  $(0\leqslant i\leqslant\inf(10,k) )$ of the below form.\\
    \begin{equation}
    \label{eq 1.1}
   \left (  \begin{array} {ccccccc}
\alpha  _{1}^{(1)} & \alpha  _{2}^{(1)}  &   \alpha_{3}^{(1)} &   \alpha_{4}^{(1)} &   \alpha_{5}^{(1)}  & \ldots  &  \alpha_{k}^{(1)} \\
\alpha  _{2}^{(1)} & \alpha  _{3}^{(1)}  &   \alpha_{4}^{(1)} &   \alpha_{5}^{(1)} &   \alpha_{6}^{(1)}  & \ldots  &  \alpha_{k+1}^{(1)} \\ 
\hline \\
\alpha  _{1}^{(2)} & \alpha  _{2}^{(2)}  &   \alpha_{3}^{(2)} &   \alpha_{4}^{(2)} &   \alpha_{5}^{(2)} & \ldots   &  \alpha_{k}^{(2)} \\
\alpha  _{2}^{(2)} & \alpha  _{3}^{(2)}  &   \alpha_{4}^{(2)} &   \alpha_{5}^{(2)}&   \alpha_{6}^{(2)}   & \ldots  &  \alpha_{k+1}^{(2)} \\ 
\hline\\
\alpha  _{1}^{(3)} & \alpha  _{2}^{(3)}  &   \alpha_{3}^{(3)}  &   \alpha_{4}^{(3)} &   \alpha_{5}^{(3)}  & \ldots  &  \alpha_{k}^{(3)} \\
\alpha  _{2}^{(3)} & \alpha  _{3}^{(3)}  &   \alpha_{4}^{(3)}&   \alpha_{5}^{(3)} &   \alpha_{6}^{(3)}  & \ldots  &  \alpha_{k+1}^{(3)} \\ 
\hline \\
\alpha  _{1}^{(4)} & \alpha  _{2}^{(4)}  &   \alpha_{3}^{(4)} &   \alpha_{4}^{(4)} &   \alpha_{5}^{(4)}  & \ldots  &  \alpha_{k}^{(4)} \\
\alpha  _{2}^{(4)} & \alpha  _{3}^{(4)}  &   \alpha_{4}^{(4)}&   \alpha_{5}^{(4)} &   \alpha_{6}^{(4)}  & \ldots  &  \alpha_{k+1}^{(4)} \\ 
\hline \\
\alpha  _{1}^{(5)} & \alpha  _{2}^{(5)}  &   \alpha_{3}^{(5)} &   \alpha_{4}^{(5)} &   \alpha_{5}^{(5)}  & \ldots  &  \alpha_{k}^{(5)} \\
\alpha  _{2}^{(5)} & \alpha  _{3}^{(5)}  &   \alpha_{4}^{(5)}&   \alpha_{5}^{(5)} &   \alpha_{6}^{(5)}  & \ldots  &  \alpha_{k+1}^{(5)}
\end{array} \right )  
\end{equation} 
We remark that this paper is based on the results in  the author's paper [12]

\section{Notation and Preliminaries}
\label{sec 2}
 \subsection{Some notations concerning the field of Laurent Series $ \mathbb{F}_{2}((T^{-1})) $ }
   We denote by $ \mathbb{F}_{2}\big(\big({T^{-1}}\big) \big)
 = \mathbb{K} $ the completion
 of the field $\mathbb{F}_{2}(T), $  the field of  rational fonctions over the
 finite field\; $\mathbb{F}_{2}$,\; for the  infinity  valuation \;
 $ \mathfrak{v}=\mathfrak{v}_{\infty }$ \;defined by \;
 $ \mathfrak{v}\big(\frac{A}{B}\big) = degB -degA $ \;
 for each pair (A,B) of non-zero polynomials.
 Then every element non-zero t in
  $\mathbb{F}_{2}\big(\big({\frac{1}{T}}\big) \big) $
 can be expanded in a unique way in a convergent Laurent series
                              $  t = \sum_{j= -\infty }^{-\mathfrak{v}(t)}t_{j}T^j
                                 \; where\; t_{j}\in \mathbb{F}_{2}. $\\
  We associate to the infinity valuation\; $\mathfrak{v}= \mathfrak{v}_{\infty }$
   the absolute value \; $\vert \cdot \vert_{\infty} $\; defined by \;
  \begin{equation*}
  \vert t \vert_{\infty} =  \vert t \vert = 2^{-\mathfrak{v}(t)}. \\
\end{equation*}
    We denote  E the  Character of the additive locally compact group
$  \mathbb{F}_{2}\big(\big({\frac{1}{T}}\big) \big) $ defined by \\
\begin{equation*}
 E\big( \sum_{j= -\infty }^{-\mathfrak{v}(t)}t_{j}T^j\big)= \begin{cases}
 1 & \text{if      }   t_{-1}= 0, \\
  -1 & \text{if      }   t_{-1}= 1.
    \end{cases}
\end{equation*}
  We denote $\mathbb{P}$ the valuation ideal in $ \mathbb{K},$ also denoted the unit interval of  $\mathbb{K},$ i.e.
  the open ball of radius 1 about 0 or, alternatively, the set of all Laurent series 
   $$ \sum_{i\geq 1}\alpha _{i}T^{-i}\quad (\alpha _{i}\in  \mathbb{F}_{2} ) $$ and, for every rational
    integer j,  we denote by $\mathbb{P}_{j} $
     the  ideal $\left\{t \in \mathbb{K}|\; \mathfrak{v}(t) > j \right\}. $
     The sets\; $ \mathbb{P}_{j}$\; are compact subgroups  of the additive
     locally compact group \; $ \mathbb{K}. $\\
      All $ t \in \mathbb{F}_{2}\Big(\Big(\frac{1}{T}\Big)\Big) $ may be written in a unique way as
$ t = [t] + \left\{t\right\}, $ \;  $  [t] \in \mathbb{F}_{2}[T] ,
 \; \left\{t\right\}\in \mathbb{P}  ( =\mathbb{P}_{0}). $\\
 We denote by dt the Haar measure on \; $ \mathbb{K} $\; chosen so that \\
  $$ \int_{\mathbb{P}}dt = 1. $$\\
  
  $$ Let \quad
  (t_{1},t_{2},\ldots,t_{n} )
 =  \big( \sum_{j=-\infty}^{-\nu(t_{1})}\alpha _{j}^{(1)}T^{j},  \sum_{j=-\infty}^{-\nu(t_{2})}\alpha _{j}^{(2)}T^{j} ,\ldots, \sum_{j=-\infty}^{-\nu(t_{n})}\alpha _{j}^{(n)}T^{j}\big) \in  \mathbb{K}^{n}. $$ 
 We denote $\psi  $  the  Character on  $(\mathbb{K}^n, +) $ defined by \\
 \begin{align*}
  \psi \big( \sum_{j=-\infty}^{-\nu(t_{1})}\alpha _{j}^{(1)}T^{j},  \sum_{j=-\infty}^{-\nu(t_{2})}\alpha _{j}^{(2)}T^{j} ,\ldots, \sum_{j=-\infty}^{-\nu(t_{n})}\alpha _{j}^{(n)}T^{j}\big) & = E \big( \sum_{j=-\infty}^{-\nu(t_{1})}\alpha _{j}^{(1)}T^{j}\big) \cdot E\big( \sum_{j=-\infty}^{-\nu(t_{2})}\alpha _{j}^{(2)}T^{j}\big)\cdots E\big(  \sum_{j=-\infty}^{-\nu(t_{n})}\alpha _{j}^{(n)}T^{j}\big) \\
  & = 
    \begin{cases}
 1 & \text{if      }     \alpha _{-1}^{(1)} +    \alpha _{-1}^{(2)}  + \ldots +   \alpha _{-1}^{(n)}   = 0 \\
  -1 & \text{if      }    \alpha _{-1}^{(1)} +    \alpha _{-1}^{(2)}  + \ldots +   \alpha _{-1}^{(n)}   =1                                                                                                                          
    \end{cases}
  \end{align*}

   \subsection{Some results concerning  n-times persymmetric matrices over  $ \mathbb{F}_{2}$}
 \label{subsec 2.2}
     $$ Set\quad
  (t_{1},t_{2},\ldots,t_{n} )
 =  \big( \sum_{i\geq 1}\alpha _{i}^{(1)}T^{-i}, \sum_{i \geq 1}\alpha  _{i}^{(2)}T^{-i},\sum_{i \geq 1}\alpha _{i}^{(3)}T^{-i},\ldots,\sum_{i \geq 1}\alpha _{i}^{(n)}T^{-i}   \big) \in  \mathbb{P}^{n}. $$

     Denote by $D^{\left[2 \atop{\vdots \atop 2}\right]\times k}(t_{1},t_{2},\ldots,t_{n} ) $
    
    the following $2n \times k $ \;  n-times  persymmetric  matrix  over the finite field  $\mathbb{F}_{2} $ 
    
  \begin{equation}
  \label{eq 2.1}
   \left (  \begin{array} {cccccccc}
\alpha  _{1}^{(1)} & \alpha  _{2}^{(1)}  &   \alpha_{3}^{(1)} &   \alpha_{4}^{(1)} &   \alpha_{5}^{(1)} &  \alpha_{6}^{(1)}  & \ldots  &  \alpha_{k}^{(1)} \\
\alpha  _{2}^{(1)} & \alpha  _{3}^{(1)}  &   \alpha_{4}^{(1)} &   \alpha_{5}^{(1)} &   \alpha_{6}^{(1)} &  \alpha_{7}^{(1)} & \ldots  &  \alpha_{k+1}^{(1)} \\ 
\hline \\
\alpha  _{1}^{(2)} & \alpha  _{2}^{(2)}  &   \alpha_{3}^{(2)} &   \alpha_{4}^{(2)} &   \alpha_{5}^{(2)} &  \alpha_{6}^{(2)} & \ldots   &  \alpha_{k}^{(2)} \\
\alpha  _{2}^{(2)} & \alpha  _{3}^{(2)}  &   \alpha_{4}^{(2)} &   \alpha_{5}^{(2)}&   \alpha_{6}^{(2)} &  \alpha_{7}^{(2)}  & \ldots  &  \alpha_{k+1}^{(2)} \\ 
\hline\\
\alpha  _{1}^{(3)} & \alpha  _{2}^{(3)}  &   \alpha_{3}^{(3)}  &   \alpha_{4}^{(3)} &   \alpha_{5}^{(3)} &  \alpha_{6}^{(3)} & \ldots  &  \alpha_{k}^{(3)} \\
\alpha  _{2}^{(3)} & \alpha  _{3}^{(3)}  &   \alpha_{4}^{(3)}&   \alpha_{5}^{(3)} &   \alpha_{6}^{(3)}  &  \alpha_{7}^{(3)} & \ldots  &  \alpha_{k+1}^{(3)} \\ 
\hline \\
\vdots & \vdots & \vdots  & \vdots  & \vdots & \vdots  & \vdots & \vdots \\
\hline \\
\alpha  _{1}^{(n)} & \alpha  _{2}^{(n)}  &   \alpha_{3}^{(n)} &   \alpha_{4}^{(n)} &   \alpha_{5}^{(n)}  &  \alpha_{6}^{(n)} & \ldots  &  \alpha_{k}^{(n)} \\
\alpha  _{2}^{(n)} & \alpha  _{3}^{(n)}  &   \alpha_{4}^{(n)}&   \alpha_{5}^{(n)} &   \alpha_{6}^{(n)}  &  \alpha_{7}^{(n)} & \ldots  &  \alpha_{k+1}^{(n)} \\ 
\end{array} \right )  
\end{equation} 
We denote by  $ \Gamma_{i}^{\left[2\atop{\vdots \atop 2}\right]\times k}$  the number of rank i n-times persymmetric matrices over $\mathbb{F}_{2}$ of the above form :  \\

  Let $ \displaystyle  f (t_{1},t_{2},\ldots,t_{n} ) $  be the exponential sum  in $ \mathbb{P}^{n} $ defined by\\
    $(t_{1},t_{2},\ldots,t_{n} ) \displaystyle\in \mathbb{P}^{n}\longrightarrow \\
    \sum_{deg Y\leq k-1}\sum_{deg U_{1}\leq  1}E(t_{1} YU_{1})
  \sum_{deg U_{2} \leq 1}E(t _{2} YU_{2}) \ldots \sum_{deg U_{n} \leq 1} E(t _{n} YU_{n}). $\vspace{0.5 cm}\\
    Then
  $$     f_{k} (t_{1},t_{2},\ldots,t_{n} ) =
  2^{2n+k- rank\big[ D^{\left[2\atop{\vdots \atop 2}\right]\times k}(t_{1},t_{2},\ldots,t_{n} )\big] } $$

    Hence  the number denoted by $ R_{q,n}^{(k)} $ of solutions \\
  
 $(Y_1,U_{1}^{(1)},U_{2}^{(1)}, \ldots,U_{n}^{(1)}, Y_2,U_{1}^{(2)},U_{2}^{(2)}, 
\ldots,U_{n}^{(2)},\ldots  Y_q,U_{1}^{(q)},U_{2}^{(q)}, \ldots,U_{n}^{(q)}   ) \in (\mathbb{F}_{2}[T])^{(n+1)q}$ \vspace{0.5 cm}\\
 of the polynomial equations  \vspace{0.5 cm}
  \[\left\{\begin{array}{c}
 Y_{1}U_{1}^{(1)} + Y_{2}U_{1}^{(2)} + \ldots  + Y_{q}U_{1}^{(q)} = 0  \\
    Y_{1}U_{2}^{(1)} + Y_{2}U_{2}^{(2)} + \ldots  + Y_{q}U_{2}^{(q)} = 0\\
    \vdots \\
   Y_{1}U_{n}^{(1)} + Y_{2}U_{n}^{(2)} + \ldots  + Y_{q}U_{n}^{(q)} = 0 
 \end{array}\right.\]
 
    $ \Leftrightarrow
    \begin{pmatrix}
   U_{1}^{(1)} & U_{1}^{(2)} & \ldots  & U_{1}^{(q)} \\ 
      U_{2}^{(1)} & U_{2}^{(2)}  & \ldots  & U_{2}^{(q)}  \\
\vdots &   \vdots & \vdots &   \vdots   \\
U_{n}^{(1)} & U_{n}^{(2)}   & \ldots  & U_{n}^{(q)} \\
 \end{pmatrix}  \begin{pmatrix}
   Y_{1} \\
   Y_{2}\\
   \vdots \\
   Y_{q} \\
  \end{pmatrix} =   \begin{pmatrix}
  0 \\
  0 \\
  \vdots \\
  0 
  \end{pmatrix} $\\

    satisfying the degree conditions \\
                   $$  degY_i \leq k-1 ,
                   \quad degU_{j}^{(i)} \leq 1, \quad  for \quad 1\leq j\leq n  \quad 1\leq i \leq q $$ \\
  is equal to the following integral over the unit interval in $ \mathbb{K}^{n} $
    $$ \int_{\mathbb{P}^{n}} f_{k}^{q}(t_{1},t_{2},\ldots,t_{n}) dt_{1}dt _{2}\ldots dt _{n}. $$
  Observing that $ f (t_{1},t_{2},\ldots,t_{n} ) $ is constant on cosets of $ \prod_{j=1}^{n}\mathbb{P}_{k+1} $ in $ \mathbb{P}^{n} $\;
  the above integral is equal to 
  
  \begin{equation}
  \label{eq 2.2}
 2^{q(2n+k) - (k+1)n}\sum_{i = 0}^{\inf(2n,k)}
  \Gamma_{i}^{\left[2\atop{\vdots \atop 2}\right]\times k} 2^{-iq} =  R_{q,n}^{(k)} 
 \end{equation}
 
 \begin{eqnarray}
 \label{eq 2.3}
\text{ Recall that $ R_{q,n}^{(k)}$ is equal to the number of solutions of the polynomial system} \nonumber \\
    \begin{pmatrix}
   U_{1}^{(1)} & U_{1}^{(2)} & \ldots  & U_{1}^{(q)} \\ 
      U_{2}^{(1)} & U_{2}^{(2)}  & \ldots  & U_{2}^{(q)}  \\
\vdots &   \vdots & \vdots &   \vdots   \\
U_{n}^{(1)} & U_{n}^{(2)}   & \ldots  & U_{n}^{(q)} \\
 \end{pmatrix}  \begin{pmatrix}
   Y_{1} \\
   Y_{2}\\
   \vdots \\
   Y_{q} \\
  \end{pmatrix} =   \begin{pmatrix}
  0 \\
  0 \\
  \vdots \\
  0 
  \end{pmatrix} \\
 \text{ satisfying the degree conditions}\nonumber \\
                     degY_i \leq k-1 ,
                   \quad degU_{j}^{(i)} \leq 1, \quad  for \quad 1\leq j\leq n  \quad 1\leq i \leq q  \nonumber
 \end{eqnarray}

 From \eqref{eq 2.2} we obtain for q = 1\\
   \begin{align}
  \label{eq 2.4}
 2^{k-(k-1)n}\sum_{i = 0}^{\inf(2n,k)}
 \Gamma_{i}^{\left[2\atop{\vdots \atop 2}\right]\times k} 2^{-i} =  R_{1,n}^{(k)} = 2^{2n}+2^k-1
  \end{align}

We have obviously \\

   \begin{align}
  \label{eq 2.5}
 \sum_{i = 0}^{k}
 \Gamma_{i}^{\left[2\atop{\vdots \atop 2}\right]\times k}  = 2^{(k+1)n}  
 \end{align}

From  the fact that the number of rank one persymmetric  matrices over $\mathbb{F}_{2}$ is equal to three  we obtain using
 combinatorial methods  : \\
 
    \begin{equation}
  \label{eq 2.6}
 \Gamma_{1}^{\left[2\atop{\vdots \atop 2}\right]\times k}  = (2^{n}-1)\cdot 3
  \end{equation}
  For more details see Cherly  [11,12]
  
  \subsection{The case n=5}

       $$ Set\quad
  (t_{1},t_{2},t_{3},t_{4},t_{5} )
 =  \big( \sum_{i\geq 1}\alpha _{i}^{(1)}T^{-i}, \sum_{i \geq 1}\alpha  _{i}^{(2)}T^{-i},\sum_{i \geq 1}\alpha _{i}^{(3)}T^{-i},\sum_{i \geq 1}\alpha _{i}^{(4)}T^{-i},\sum_{i \geq 1}\alpha _{i}^{(5)}T^{-i}   \big) \in  \mathbb{P}^{5}. $$

     Denote by $D^{\left[2 \atop{2\atop{2 \atop {2\atop 2}}}\right]\times k}(t_{1},t_{2},t_{3},t_{4},t_{5}) $
    
    the following $10 \times k $ \; quintuple  persymmetric  matrix  over the finite field  $\mathbb{F}_{2} $ 
   
  \begin{displaymath}
   \left (  \begin{array} {cccccccc}
\alpha  _{1}^{(1)} & \alpha  _{2}^{(1)}  &   \alpha_{3}^{(1)} &   \alpha_{4}^{(1)} &   \alpha_{5}^{(1)}   & \ldots  &  \alpha_{k}^{(1)} \\
\alpha  _{2}^{(1)} & \alpha  _{3}^{(1)}  &   \alpha_{4}^{(1)} &   \alpha_{5}^{(1)} &   \alpha_{6}^{(1)} & \ldots  &  \alpha_{k+1}^{(1)} \\ 
\hline \\
\alpha  _{1}^{(2)} & \alpha  _{2}^{(2)}  &   \alpha_{3}^{(2)} &   \alpha_{4}^{(2)} &   \alpha_{5}^{(2)} & \ldots   &  \alpha_{k}^{(2)} \\
\alpha  _{2}^{(2)} & \alpha  _{3}^{(2)}  &   \alpha_{4}^{(2)} &   \alpha_{5}^{(2)}&   \alpha_{6}^{(2)}   & \ldots  &  \alpha_{k+1}^{(2)} \\ 
\hline\\
\alpha  _{1}^{(3)} & \alpha  _{2}^{(3)}  &   \alpha_{3}^{(3)}  &   \alpha_{4}^{(3)} &   \alpha_{5}^{(3)} & \ldots  &  \alpha_{k}^{(3)} \\
\alpha  _{2}^{(3)} & \alpha  _{3}^{(3)}  &   \alpha_{4}^{(3)}&   \alpha_{5}^{(3)} &   \alpha_{6}^{(3)} & \ldots  &  \alpha_{k+1}^{(3)} \\ 
\hline \\
\alpha  _{1}^{(4)} & \alpha  _{2}^{(4)}  &   \alpha_{3}^{(4)} &   \alpha_{4}^{(4)} &   \alpha_{5}^{(4)}  & \ldots  &  \alpha_{k}^{(4)} \\
\alpha  _{2}^{(4)} & \alpha  _{3}^{(4)}  &   \alpha_{4}^{(4)}&   \alpha_{5}^{(4)} &   \alpha_{6}^{(4)}  & \ldots  &  \alpha_{k+1}^{(4)} \\ 
\hline \\
\alpha  _{1}^{(5)} & \alpha  _{2}^{(5)}  &   \alpha_{3}^{(5)} &   \alpha_{4}^{(5)} &   \alpha_{5}^{(5)} & \ldots  &  \alpha_{k}^{(4)} \\
\alpha  _{2}^{(5)} & \alpha  _{3}^{(5)}  &   \alpha_{4}^{(5)}&   \alpha_{5}^{(5)} &   \alpha_{6}^{(5)}  & \ldots  &  \alpha_{k+1}^{(4)} 
\end{array} \right )  
\end{displaymath} 
We denote by  $ \Gamma_{i}^{\left[2\atop{2 \atop{ 2\atop {2\atop 2}}}\right]\times k}$  the number of rank i quintuple persymmetric matrices over $\mathbb{F}_{2}$ of the above form :  \\

  Let $ \displaystyle  f (t_{1},t_{2},t_{3},t_{4},t_{5} ) $  be the exponential sum  in $ \mathbb{P}^{5} $ defined by\\
    $(t_{1},t_{2},t_{3},t_{4},t_{5}) \displaystyle\in \mathbb{P}^{5}\longrightarrow \\
    \sum_{deg Y\leq k-1}\sum_{deg U_{1}\leq  1}E(t_{1} YU_{1})
  \sum_{deg U_{2} \leq 1}E(t _{2} YU_{2}) \sum_{deg U_{3} \leq 1}E(t _{3} YU_{3}) \sum_{deg U_{4} \leq 1} E(t _{4} YU_{4}) \sum_{deg U_{5} \leq 1} E(t _{5} YU_{5}). $\vspace{0.5 cm}\\
    Then
  $$     f_{k} (t_{1},t_{2},t_{3},t_{4},t_{5} ) =
  2^{10+k- rank\big[ D^{\left[2\atop{2\atop{2 \atop {2\atop 2}}}\right]\times k}(t_{1},t_{2},t_{3},t_{4},t_{5} )\big] } $$

    Hence  the number denoted by $ R_{q,5}^{(k)} $ of solutions \\
  
 $(Y_1,U_{1}^{(1)},U_{2}^{(1)},U_{3}^{(1)} ,U_{4}^{(1)}, U_{5}^{(1)},Y_2,U_{1}^{(2)},U_{2}^{(2)}, 
U_{3}^{(2)},U_{4}^{(2)},U_{5}^{(2)},\ldots  Y_q,U_{1}^{(q)},U_{2}^{(q)}, U_{3}^{(q)},U_{4}^{(q)},U_{5}^{(q)}   ) \in (\mathbb{F}_{2}[T])^{6q}$ \vspace{0.5 cm}\\
 of the polynomial equations  \vspace{0.5 cm}
  \[\left\{\begin{array}{c}
 Y_{1}U_{1}^{(1)} + Y_{2}U_{1}^{(2)} + \ldots  + Y_{q}U_{1}^{(q)} = 0  \\
    Y_{1}U_{2}^{(1)} + Y_{2}U_{2}^{(2)} + \ldots  + Y_{q}U_{2}^{(q)} = 0\\
    Y_{1}U_{3}^{(1)} + Y_{3}U_{3}^{(2)} + \ldots  + Y_{q}U_{3}^{(q)} = 0\\ 
   Y_{1}U_{4}^{(1)} + Y_{2}U_{4}^{(2)} + \ldots  + Y_{q}U_{4}^{(q)} = 0 \\
  Y_{1}U_{5}^{(1)} + Y_{2}U_{5}^{(2)} + \ldots  + Y_{q}U_{5}^{(q)} = 0  
 \end{array}\right.\]
 
    $ \Leftrightarrow
    \begin{pmatrix}
   U_{1}^{(1)} & U_{1}^{(2)} & \ldots  & U_{1}^{(q)} \\ 
      U_{2}^{(1)} & U_{2}^{(2)}  & \ldots  & U_{2}^{(q)}  \\
 U_{3}^{(1)} & U_{3}^{(2)}  & \ldots  & U_{3}^{(q)}  \\
U_{4}^{(1)} & U_{4}^{(2)}   & \ldots  & U_{4}^{(q)} \\
U_{5}^{(1)} & U_{5}^{(2)}   & \ldots  & U_{5}^{(q)}
 \end{pmatrix}  \begin{pmatrix}
   Y_{1} \\
   Y_{2}\\
   \vdots \\
   Y_{q} \\
  \end{pmatrix} =   \begin{pmatrix}
  0 \\
  0 \\
  0 \\
  0 \\
  0
  \end{pmatrix} $\\

    satisfying the degree conditions \\
                   $$  degY_i \leq k-1 ,
                   \quad degU_{j}^{(i)} \leq 1, \quad  for \quad 1\leq j\leq 5  \quad 1\leq i \leq q $$ \\
  is equal to the following integral over the unit interval in $ \mathbb{K}^{5} $
    $$ \int_{\mathbb{P}^{5}} f_{k}^{q}(t_{1},t_{2},t_{3},t_{4},t_{5}) dt_{1}dt _{2}dt_{3} dt _{4}dt_{5}$$
  Observing that $ f (t_{1},t_{2},t_{3},t_{4},t_{5} ) $ is constant on cosets of $ \prod_{j=1}^{5}\mathbb{P}_{k+1} $ in $ \mathbb{P}^{5} $\;
  the above integral is equal to 
  
  \begin{equation}
  \label{eq 2.7}
 2^{q(10+k) - 5(k+1)}\sum_{i = 0}^{\inf{(10,k)}}
  \Gamma_{i}^{\left[2\atop{2\atop{2 \atop {2\atop 2}}}\right]\times k} 2^{-iq} =  R_{q,5}^{(k)} \quad \text{where} \; k\geqslant 1
 \end{equation}
 \newpage
 
We shall need the following results.\\

 \textbf{Result 1} : 
  \newline
 We have whenever $k\geqslant 4:$ See Cherly [12]  \\
    \begin{equation}
    \label{eq 2.8}
 \Gamma_{i}^{\left[2\atop {2\atop {2\atop2 }}\right]\times k}  = \begin{cases}
1 & \text{if  } i = 0,        \\
 45 & \text{if   } i=1,\\
 30\cdot2^{k}+1410  & \text{if   }  i = 2,  \\
 1470\cdot 2^{k}+ 31920  & \text{if   }  i = 3, \\
  140\cdot 2^{2k} +42420\cdot 2^{k}+276640  & \text{if   }  i=4, \\  
  6300 \cdot 2^{2k}+630000 \cdot2^{k} -11692800 & \text{if   }  i=5, \\
 120\cdot 2^{3k}+123480\cdot2^{2k} -6142080\cdot 2^{k}  +66170880    & \text{if   }  i=6. \\
 3720\cdot 2^{3k}-416640\cdot2^{2k}+13332480\cdot2^{k} -121896960 & \text{if   }  i=7. \\
    16\cdot 2^{4k}- 3840\cdot 2^{3k} +286720\cdot2^{2k}-7864320\cdot 2^{k} +2^{26}& \text{if   }  i=8.    
  \end{cases}    
  \end{equation}
  where    $ \Gamma_{i}^{\left[2\atop{2 \atop{ 2\atop 2}}\right]\times k}$ denotes the number 
  of  quadruple  persymmetric matrices in   $ \mathbb{F}_{2}$ of rank i  $(0\leqslant i\leqslant\inf(8,k) )$ of the below form.\\
    \begin{equation*}
    \left (  \begin{array} {ccccccc}
\alpha  _{1}^{(1)} & \alpha  _{2}^{(1)}  &   \alpha_{3}^{(1)} &   \alpha_{4}^{(1)} &   \alpha_{5}^{(1)} & \ldots  &  \alpha_{k}^{(1)} \\
\alpha  _{2}^{(1)} & \alpha  _{3}^{(1)}  &   \alpha_{4}^{(1)} &   \alpha_{5}^{(1)} &   \alpha_{6}^{(1)}  & \ldots  &  \alpha_{k+1}^{(1)} \\ 
\hline \\
\alpha  _{1}^{(2)} & \alpha  _{2}^{(2)}  &   \alpha_{3}^{(2)} &   \alpha_{4}^{(2)} &   \alpha_{5}^{(2)}& \ldots   &  \alpha_{k}^{(2)} \\
\alpha  _{2}^{(2)} & \alpha  _{3}^{(2)}  &   \alpha_{4}^{(2)} &   \alpha_{5}^{(2)}&   \alpha_{6}^{(2)}  & \ldots  &  \alpha_{k+1}^{(2)} \\ 
\hline\\
\alpha  _{1}^{(3)} & \alpha  _{2}^{(3)}  &   \alpha_{3}^{(3)}  &   \alpha_{4}^{(3)} &   \alpha_{5}^{(3)}  & \ldots  &  \alpha_{k}^{(3)} \\
\alpha  _{2}^{(3)} & \alpha  _{3}^{(3)}  &   \alpha_{4}^{(3)}&   \alpha_{5}^{(3)} &   \alpha_{6}^{(3)} & \ldots  &  \alpha_{k+1}^{(3)} \\ 
\hline \\
\alpha  _{1}^{(4)} & \alpha  _{2}^{(4)}  &   \alpha_{3}^{(4)} &   \alpha_{4}^{(4)} &   \alpha_{5}^{(4)} & \ldots  &  \alpha_{k}^{(4)} \\
\alpha  _{2}^{(4)} & \alpha  _{3}^{(4)}  &   \alpha_{4}^{(4)}&   \alpha_{5}^{(4)} &   \alpha_{6}^{(4)}  & \ldots  &  \alpha_{k+1}^{(4)} \\ 
\end{array} \right )  
\end{equation*} 
\textbf{Result 2}
  \newline
  
 The  $ \Gamma_{i}^{\left[2\atop{\vdots \atop 2}\right]\times k} $  where  $ 0\leq i\leq \inf(2n,k)$ (see subsection  \ref{subsec 2.2} )
are solutions to the below  system.  See Cherly[12 ]\\

  \begin{equation}
  \label{eq 2.9}
 \begin{cases} 
 \displaystyle  \Gamma_{0}^{\left[2\atop{\vdots \atop 2}\right]\times k}  = 1 \quad \text{if} \quad  k\geqslant 1 \\
\displaystyle  \Gamma_{1}^{\left[2\atop{\vdots \atop 2}\right]\times k}  = (2^{n}-1)\cdot 3 \quad \text{if} \quad  k\geqslant 2 \\
\displaystyle \Gamma_{2}^{\left[2\atop{\vdots \atop 2}\right]\times k} = 7\cdot2^{2n}+(2^{k+1}-25) \cdot 2^{n}-2^{k+1}+18 \quad \text{for} \quad k\geqslant 3\\
\displaystyle  \Gamma_{3}^{\left[2\atop{\vdots \atop 2}\right]\times k} = 15\cdot2^{3n} + (7\cdot2^k-133)\cdot2^{2n}+ (294-21\cdot 2^k) \cdot 2^{n}   -176+14\cdot2^k \quad \text{for} \quad k\geqslant 4\\
\displaystyle  \Gamma_{4}^{\left[2\atop{\vdots \atop 2}\right]\times k} = 31\cdot2^{4n} + \frac{35\cdot2^{k}-1210}{2}\cdot2^{3n}
+ \frac{2^{2k+2}-783\cdot2^{k}+19028}{6}\cdot 2^{2n}\\
\displaystyle +(-2^{2k+1}+269\cdot2^{k}-5744)\cdot2^n 
 +\frac{2^{2k+2}-117\cdot2^{k+2}+9440}{3}  \quad \text{for} \quad k\geqslant 5\\
   \displaystyle  \Gamma_{5}^{\left[2\atop{\vdots \atop 2}\right]\times k} = 63\cdot2^{5n} + (\frac{155}{4}\cdot2^{k}-2573)\cdot2^{4n}
+ (\frac{5}{2}\cdot2^{2k}-\frac{2565}{4}\cdot2^{k}+29150)\cdot2^{3n}\\
\displaystyle +\frac{1}{2}\cdot(-35\cdot2^{2k}+6265\cdot2^{k}-247520)\cdot2^{2n} 
\displaystyle +(35\cdot2^{2k}-5490\cdot2^{k}+203872)\cdot2^{n}\\
-20\cdot2^{2k}+2960\cdot2^{k}-106752  \quad \text{for} \quad k\geqslant 6\\
\displaystyle  \sum_{i = 0}^{\inf(2n,k)} \Gamma_{i}^{\left[2\atop{\vdots \atop 2}\right]\times k}  = 2^{(k+1)n} \\ 
  \displaystyle  \sum_{i = 0}^{\inf(2n,k)} \Gamma_{i}^{\left[2\atop{\vdots \atop 2}\right]\times k} 2^{-i}  = 2^{n+k(n-1)}+2^{(k-1)n}-2^{(k-1)n-k}\\
  \displaystyle \sum_{i = 0}^{\inf(2n,k)} \Gamma_{i}^{\left[2\atop{\vdots \atop 2}\right]\times k} 2^{-2i}  =
   2^{n+k(n-2)}+2^{-n+k(n-2)}\cdot[3\cdot2^k-3] +2^{-2n+k(n-2)}\cdot[6\cdot2^{k-1}-6] \\
   +2^{-3n+kn}-6\cdot2^{n(k-3)-k}+8\cdot2^{-3n+k(n-2)}
\end{cases}
    \end{equation}
  \textbf{Result 3}
    \newline
The number of rank 10 quintuple persymmetric matrices
   of the  form \eqref{eq 1.1} is equal to :\\
  $ \displaystyle 2^{5}\prod_{j=1}^{5}(2^{k}-2^{10 -j}) .$See Cherly[10, section 2 ] \\
  That is :\\
  \begin{equation}
\label{eq 2.10}
 \displaystyle  \Gamma_{10}^{\left[2\atop {2\atop {2\atop{2\atop 2} }}\right]\times k} = 2^{5}\prod_{j=1}^{5}(2^{k}-2^{10 -j}) 
\end{equation}

  \subsection{Computation of the number of quintuple persymmetric matrices of the form (1.1) of rank I}
\begin{thm}
\label{thm 2.1}
We have whenever $k\geqslant 5: $\\
 \begin{equation}
 \label{eq 2.11}
  \begin{cases} 
 \displaystyle   \Gamma_{0}^{\left[2\atop {2\atop {2\atop{2\atop2 }}}\right]\times k}   = 1 \quad \text{if} \quad  k\geqslant 1 \\
\displaystyle    \Gamma_{1}^{\left[2\atop {2\atop {2\atop{2\atop2 }}}\right]\times k}   = 93 \quad \text{if} \quad  k\geqslant 2 \\
\displaystyle   \Gamma_{2}^{\left[2\atop {2\atop {2\atop{2\atop2 }}}\right]\times k}  = 31\cdot2^{k+1}+6386 \quad \text{for} \quad k\geqslant 3\\
\displaystyle    \Gamma_{3}^{\left[2\atop {2\atop {2\atop{2\atop2 }}}\right]\times k}  = 6510\cdot2^{k}+364560 \quad \text{for} \quad k\geqslant 4\\
\displaystyle    \Gamma_{4}^{\left[2\atop {2\atop {2\atop{2\atop2 }}}\right]\times k}  = 620\cdot2^{2k}+448260\cdot2^{k}+15748000 \quad \text{for} \quad k\geqslant 5\\
  \displaystyle    \Gamma_{5}^{\left[2\atop {2\atop {2\atop{2\atop2 }}}\right]\times k}  =  65100\cdot2^{2k}+22654800\cdot2^{k}+250817280   \quad \text{for} \quad k\geqslant 6\\
   \displaystyle   \Gamma_{6}^{\left[2\atop {2\atop {2\atop{2\atop2 }}}\right]\times k}  = 1240\cdot[2^{3k}
  +3199\cdot2^{2k}+2^{7}\cdot2913\cdot2^{k}-18883\cdot 2^{10}] \quad \text{for} \quad k\geqslant 7 \\
   \displaystyle   \Gamma_{7}^{\left[2\atop {2\atop {2\atop{2\atop2 }}}\right]\times k}  = 115320\cdot[2^{3k}
  +1148\cdot2^{2k}-2^7\cdot 917\cdot2^{k}+311\cdot2^{13}] \quad \text{for} \quad k\geqslant 8 \\
  \displaystyle   \Gamma_{8}^{\left[2\atop {2\atop {2\atop{2\atop2 }}}\right]\times k}  = 496\cdot[2^{4k}+9525\cdot2^{3k}
     -2169440\cdot2^{2k}+68115\cdot2^{11}\cdot2^{k}-9749\cdot2^{18}] \quad \text{for} \quad k\geqslant 9 \\
   \displaystyle   \Gamma_{9}^{\left[2\atop {2\atop {2\atop{2\atop2 }}}\right]\times k}  = 31248\cdot[2^{4k}-480\cdot2^{3k}
  +71680\cdot2^{2k}-3932160\cdot2^{k}+2^{26}]  \quad \text{for} \quad k\geqslant 10 \\
 \displaystyle   \Gamma_{10}^{\left[2\atop {2\atop {2\atop{2\atop2 }}}\right]\times k}  = 2^5\cdot[2^{5k}-992\cdot2^{4k}+317440\cdot2^{3k}-40632320\cdot2^{2k}+2080374784\cdot2^{k}-2^{35}]  \; \text{for} \; k\geqslant 10 \\
\end{cases}
    \end{equation}
\end{thm}
\newpage
\begin{proof}
\textbf{Step 1}
\newline
From the expressions of  $ \Gamma_{i}^{\left[2\atop {2\atop {2\atop2 }}\right]\times k} $ in \eqref{eq 2.8} we postulate that 
   \begin{equation}
    \label{eq 2.12}
 \Gamma_{i}^{\left[2\atop {2\atop {2\atop{2\atop 2}}}\right]\times k}  = \begin{cases}
1 & \text{if  } i = 0,        \\
 a_{1} & \text{if   } i=1,\\  
 a_{2}\cdot2^{k}+ b_{2}  & \text{if   }  i = 2,  \\
 a_{3}\cdot 2^{k}+ b_{3}  & \text{if   }  i = 3, \\
  a_{4}\cdot 2^{2k} +b_{4}\cdot 2^{k}+c_{4}  & \text{if   }  i=4, \\  
  a_{5} \cdot 2^{2k}+b_{5} \cdot2^{k} +c_{5} & \text{if   }  i=5, \\
  a_{6}\cdot 2^{3k}+b_{6}\cdot2^{2k} +c_{6}\cdot 2^{k}  + d_{6}   & \text{if   }  i=6. \\
a_{7}\cdot 2^{3k}+b_{7}\cdot2^{2k}+c_{7}\cdot2^{k} +d_{7} & \text{if   }  i=7. \\
 a_{8}\cdot 2^{4k}+b_{8}\cdot 2^{3k} +c_{8}\cdot2^{2k}+d_{8}\cdot 2^{k} +e_{8} & \text{if   }  i=8.\\
     a_{9}\cdot 2^{4k}+b_{9}\cdot 2^{3k} +c_{9}\cdot2^{2k}+d_{9}\cdot 2^{k} +e_{9} & \text{if   }  i=9.\\   
a_{10}\cdot2^{5k} + b_{10}\cdot 2^{4k}+c_{10}\cdot 2^{3k} +d_{10}\cdot2^{2k}+e_{10}\cdot 2^{k} +f_{10} & \text{if   }  i=10.
  \end{cases}    
  \end{equation}
\textbf{Step 2}

	Equally we postulate that :\\
 \begin{equation}
 \label{eq 2.13}
  \begin{cases} 
    \displaystyle   \Gamma_{6}^{\left[2\atop {2\atop {2\atop{2\atop2 }}}\right]\times k}  = 0 \quad \text{for} \quad k = 5 \\
   \displaystyle   \Gamma_{7}^{\left[2\atop {2\atop {2\atop{2\atop2 }}}\right]\times k}  = 0\quad \text{for} \quad k\in \{5,6\} \\
   \displaystyle   \Gamma_{8}^{\left[2\atop {2\atop {2\atop{2\atop2 }}}\right]\times k} =0\quad \text{for} \quad k\in \{5,6,7\} \\
 \displaystyle   \Gamma_{9}^{\left[2\atop {2\atop {2\atop{2\atop2 }}}\right]\times k} = 0 \quad \text{for} \quad k\in \{5,6,7,8\} \\
   \displaystyle   \Gamma_{10}^{\left[2\atop {2\atop {2\atop{2\atop2 }}}\right]\times k}= 0 \quad \text{for} \quad k\in \{5,6,7,8,9\} 
\end{cases}
    \end{equation}
\textbf{Step 3}
Combining  \eqref{eq 2.9} with n=5 , \eqref{eq 2.10} and \eqref{eq 2.12} we obtain :\\
 \begin{equation}
 \label{eq 2.14}
  \begin{cases} 
 \displaystyle   \Gamma_{0}^{\left[2\atop {2\atop {2\atop{2\atop2 }}}\right]\times k}   = 1 \quad \text{for} \quad  k\geqslant 1 \\
\displaystyle    \Gamma_{1}^{\left[2\atop {2\atop {2\atop{2\atop2 }}}\right]\times k}   = 93 \quad \text{for} \quad  k\geqslant 2 \\
\displaystyle   \Gamma_{2}^{\left[2\atop {2\atop {2\atop{2\atop2 }}}\right]\times k}  = 31\cdot2^{k+1}+6386 \quad \text{for} \quad k\geqslant 3\\
\displaystyle    \Gamma_{3}^{\left[2\atop {2\atop {2\atop{2\atop2 }}}\right]\times k}  = 6510\cdot2^{k}+364560 \quad \text{for} \quad k\geqslant 4\\
\displaystyle    \Gamma_{4}^{\left[2\atop {2\atop {2\atop{2\atop2 }}}\right]\times k}  = 620\cdot2^{2k}+448260\cdot2^{k}+15748000 \quad \text{for} \quad k\geqslant 5\\
  \displaystyle    \Gamma_{5}^{\left[2\atop {2\atop {2\atop{2\atop2 }}}\right]\times k}  =  65100\cdot2^{2k}+22654800\cdot2^{k}+250817280   \quad \text{for} \quad k\geqslant 6\\
   \displaystyle    \Gamma_{6}^{\left[2\atop {2\atop {2\atop{2\atop2 }}}\right]\times k}  =  a_{6}\cdot 2^{3k}+b_{6}\cdot2^{2k} +c_{6}\cdot 2^{k}  + d_{6}   \quad \text{for  }   \quad  k\geqslant 7 \\
 \displaystyle    \Gamma_{7}^{\left[2\atop {2\atop {2\atop{2\atop2 }}}\right]\times k} = a_{7}\cdot 2^{3k}+b_{7}\cdot2^{2k}+c_{7}\cdot2^{k} +d_{7}  \quad \text{for  } \quad  k\geqslant 8  \\
 \displaystyle    \Gamma_{8}^{\left[2\atop {2\atop {2\atop{2\atop2 }}}\right]\times k} = a_{8}\cdot 2^{4k}+b_{8}\cdot 2^{3k} +c_{8}\cdot2^{2k}+d_{8}\cdot 2^{k} +e_{8} \quad  \text{for   } \quad  k\geqslant 9  \\
   \displaystyle    \Gamma_{9}^{\left[2\atop {2\atop {2\atop{2\atop2 }}}\right]\times k}  =  a_{9}\cdot 2^{4k}+b_{9}\cdot 2^{3k} +c_{9}\cdot2^{2k}+d_{9}\cdot 2^{k} +e_{9} \quad  \text{for   }\quad  k\geqslant 10   \\   
 \displaystyle   \Gamma_{10}^{\left[2\atop {2\atop {2\atop{2\atop2 }}}\right]\times k}  = 2^5\cdot[2^{5k}-992\cdot2^{4k}+317440\cdot2^{3k}-40632320\cdot2^{2k}+2080374784\cdot2^{k}-2^{35}]  \; \text{for} \; k\geqslant 10 \\
\end{cases}
    \end{equation}
   and the relations:\\
  \begin{equation}
 \label{eq 2.15}
  \begin{cases} 
 \displaystyle  \sum_{i = 0}^{10} \Gamma_{i}^{\left[2\atop {2\atop {2\atop{2\atop2 }}}\right]\times k} = 2^{5k+5}\\
   \displaystyle  \sum_{i = 0}^{10} \Gamma_{i}^{\left[2\atop {2\atop {2\atop{2\atop2 }}}\right]\times k}  2^{10-i}  = 2^{5k+5}+1023\cdot2^{4k+5}\\  
  \displaystyle \sum_{i = 0}^{10}  \Gamma_{i}^{\left[2\atop {2\atop {2\atop{2\atop2 }}}\right]\times k}  2^{20-2i}  =
   2^{5k+5}+3162\cdot2^{4k+5}+1045320\cdot2^{3k+5}
\end{cases}
    \end{equation}
   
\textbf{Step 4}
\newline

Computation of $a_{8},a_{9}$ in \eqref{eq 2.14}.\\
From \eqref{eq 2.14} and\eqref{eq 2.15} we get:\\
 \begin{equation*}
   \displaystyle 
    \begin{pmatrix}
   1 &1  \\ 
    4 & 2   \\
   16 &  4 \\
   \end{pmatrix}\displaystyle  \begin{pmatrix}
   a_{8} \\
   a_{9} \\
 \end{pmatrix} =   \begin{pmatrix}
  992\cdot 2^5 \\
  992\cdot2^5+1023\cdot2^5\\
 992\cdot2^5+3162\cdot2^5  
   \end{pmatrix} 
    \end{equation*}
    \begin{equation} 
    \label{eq 2.16}
\displaystyle  \Leftrightarrow
    \begin{pmatrix}
   a_{8} \\
   a_{9} \\
    \end{pmatrix} = \begin{pmatrix}
 496\\
 31248
  \end{pmatrix} 
  \end{equation}

\textbf{Step 5}
\newline

Computation of  $ \Gamma_{9}^{\left[2\atop {2\atop {2\atop{2\atop2 }}}\right]\times k} $ in \eqref{eq 2.14}.\\
From \eqref{eq 2.13} and \eqref{eq 2.16} we obtain :\\

\begin{equation}
\label{eq 2.17}
 \Gamma_{9}^{\left[2\atop {2\atop {2\atop{2\atop2 }}}\right]\times k} = a_{9}(2^k-2^5)(2^k-2^6)(2^k-2^7)(2^k-2^8)
=31248\cdot[2^{4k}-480\cdot2^{3k} +71680\cdot2^{2k}-3932160\cdot2^{k}+2^{26}] 
\end{equation}
To sum up we deduce from \eqref{eq 2.17},\eqref{eq 2.16} and \eqref{eq 2.14} \\
 \begin{equation}
 \label{eq 2.18}
  \begin{cases} 
 \displaystyle   \Gamma_{0}^{\left[2\atop {2\atop {2\atop{2\atop2 }}}\right]\times k}   = 1 \quad \text{for} \quad  k\geqslant 1 \\
\displaystyle    \Gamma_{1}^{\left[2\atop {2\atop {2\atop{2\atop2 }}}\right]\times k}   = 93 \quad \text{for} \quad  k\geqslant 2 \\
\displaystyle   \Gamma_{2}^{\left[2\atop {2\atop {2\atop{2\atop2 }}}\right]\times k}  = 31\cdot2^{k+1}+6386 \quad \text{for} \quad k\geqslant 3\\
\displaystyle    \Gamma_{3}^{\left[2\atop {2\atop {2\atop{2\atop2 }}}\right]\times k}  = 6510\cdot2^{k}+364560 \quad \text{for} \quad k\geqslant 4\\
\displaystyle    \Gamma_{4}^{\left[2\atop {2\atop {2\atop{2\atop2 }}}\right]\times k}  = 620\cdot2^{2k}+448260\cdot2^{k}+15748000 \quad \text{for} \quad k\geqslant 5\\
  \displaystyle    \Gamma_{5}^{\left[2\atop {2\atop {2\atop{2\atop2 }}}\right]\times k}  =  65100\cdot2^{2k}+22654800\cdot2^{k}+250817280   \quad \text{for} \quad k\geqslant 6\\
   \displaystyle    \Gamma_{6}^{\left[2\atop {2\atop {2\atop{2\atop2 }}}\right]\times k}  =  a_{6}\cdot 2^{3k}+b_{6}\cdot2^{2k} +c_{6}\cdot 2^{k}  + d_{6}   \quad \text{for  }   \quad  k\geqslant 7 \\
 \displaystyle    \Gamma_{7}^{\left[2\atop {2\atop {2\atop{2\atop2 }}}\right]\times k} = a_{7}\cdot 2^{3k}+b_{7}\cdot2^{2k}+c_{7}\cdot2^{k} +d_{7}  \quad \text{for  } \quad  k\geqslant 8  \\
 \displaystyle    \Gamma_{8}^{\left[2\atop {2\atop {2\atop{2\atop2 }}}\right]\times k} = 496\cdot 2^{4k}+b_{8}\cdot 2^{3k} +c_{8}\cdot2^{2k}+d_{8}\cdot 2^{k} +e_{8} \quad  \text{for   } \quad  k\geqslant 9  \\
 \displaystyle     \Gamma_{9}^{\left[2\atop {2\atop {2\atop{2\atop2 }}}\right]\times k} =31248\cdot[2^{4k}-480\cdot2^{3k} +71680\cdot2^{2k}-3932160\cdot2^{k}+2^{26}]  \quad  \text{for   }\quad  k\geqslant 10   \\   
 \displaystyle   \Gamma_{10}^{\left[2\atop {2\atop {2\atop{2\atop2 }}}\right]\times k}  = 2^5\cdot[2^{5k}-992\cdot2^{4k}+317440\cdot2^{3k}-40632320\cdot2^{2k}+2080374784\cdot2^{k}-2^{35}]  \; \text{for} \; k\geqslant 10 \\
\end{cases}
    \end{equation}
\textbf{Step 6}
\newline
Computation of $a_{6},a_{7} \; \text{and}\; b_{8}$ in \eqref{eq 2.18}.\\
From \eqref{eq 2.18} and \eqref{eq 2.15} we get:\\
 \begin{equation*}
   \displaystyle 
    \begin{pmatrix}
   1 &1 & 1 \\ 
    2^4 & 2^3 &2^2   \\
   2^8 &  2^6 & 2^4 \\
   \end{pmatrix}\displaystyle  \begin{pmatrix}
   a_{6} \\
   a_{7} \\
   b_{8}
 \end{pmatrix} =   \begin{pmatrix}
  31248\cdot 480 -317440\cdot2^5 \\
 2\cdot 31248\cdot 480 -2^5\cdot317440\\
  2^2\cdot 31248\cdot 480 -2^5\cdot 317440 +1045320\cdot2^5 \\
   \end{pmatrix} 
    \end{equation*}
    \begin{equation} 
    \label{eq 2.19}
\displaystyle  \Leftrightarrow
    \begin{pmatrix}
   a_{6} \\
   a_{7} \\
   b_{8}
    \end{pmatrix} = \begin{pmatrix}
 1240\\
 115320\\
 496\cdot 9525
  \end{pmatrix} 
  \end{equation}
To sum up we deduce from \eqref{eq 2.19},\eqref{eq 2.18} \\

 \begin{equation}
 \label{eq 2.20}
  \begin{cases} 
 \displaystyle   \Gamma_{0}^{\left[2\atop {2\atop {2\atop{2\atop2 }}}\right]\times k}   = 1 \quad \text{for} \quad  k\geqslant 1 \\
\displaystyle    \Gamma_{1}^{\left[2\atop {2\atop {2\atop{2\atop2 }}}\right]\times k}   = 93 \quad \text{for} \quad  k\geqslant 2 \\
\displaystyle   \Gamma_{2}^{\left[2\atop {2\atop {2\atop{2\atop2 }}}\right]\times k}  = 31\cdot2^{k+1}+6386 \quad \text{for} \quad k\geqslant 3\\
\displaystyle    \Gamma_{3}^{\left[2\atop {2\atop {2\atop{2\atop2 }}}\right]\times k}  = 6510\cdot2^{k}+364560 \quad \text{for} \quad k\geqslant 4\\
\displaystyle    \Gamma_{4}^{\left[2\atop {2\atop {2\atop{2\atop2 }}}\right]\times k}  = 620\cdot2^{2k}+448260\cdot2^{k}+15748000 \quad \text{for} \quad k\geqslant 5\\
  \displaystyle    \Gamma_{5}^{\left[2\atop {2\atop {2\atop{2\atop2 }}}\right]\times k}  =  65100\cdot2^{2k}+22654800\cdot2^{k}+250817280   \quad \text{for} \quad k\geqslant 6\\
   \displaystyle    \Gamma_{6}^{\left[2\atop {2\atop {2\atop{2\atop2 }}}\right]\times k}  =  1240\cdot 2^{3k}+b_{6}\cdot2^{2k} +c_{6}\cdot 2^{k}  + d_{6}   \quad \text{for  }   \quad  k\geqslant 7 \\
 \displaystyle    \Gamma_{7}^{\left[2\atop {2\atop {2\atop{2\atop2 }}}\right]\times k} = 115320\cdot 2^{3k}+b_{7}\cdot2^{2k}+c_{7}\cdot2^{k} +d_{7}  \quad \text{for  } \quad  k\geqslant 8  \\
 \displaystyle    \Gamma_{8}^{\left[2\atop {2\atop {2\atop{2\atop2 }}}\right]\times k} = 496\cdot 2^{4k}+496\cdot9525\cdot 2^{3k} +c_{8}\cdot2^{2k}+d_{8}\cdot 2^{k} +e_{8} \quad  \text{for   } \quad  k\geqslant 9  \\
 \displaystyle     \Gamma_{9}^{\left[2\atop {2\atop {2\atop{2\atop2 }}}\right]\times k} =31248\cdot[2^{4k}-480\cdot2^{3k} +71680\cdot2^{2k}-3932160\cdot2^{k}+2^{26}]  \quad  \text{for   }\quad  k\geqslant 10   \\   
 \displaystyle   \Gamma_{10}^{\left[2\atop {2\atop {2\atop{2\atop2 }}}\right]\times k}  = 2^5\cdot[2^{5k}-992\cdot2^{4k}+317440\cdot2^{3k}-40632320\cdot2^{2k}+2080374784\cdot2^{k}-2^{35}]  \; \text{for} \; k\geqslant 10 \\
\end{cases}
    \end{equation}

\textbf{Step 7}
\newline
Computation of $c_{8},d_{8} \; \text{and}\; e_{8}$ in \eqref{eq 2.20}.\\
From \eqref{eq 2.13} we get :  \\
 \begin{equation*}
   \displaystyle 
    \begin{pmatrix}
   2^{10} & 2^5 & 1 \\ 
    2^{12} & 2^{6} & 1   \\
   2^{14} &  2^7 & 1 \\
   \end{pmatrix}\displaystyle  \begin{pmatrix}
   c_{8} \\
   d_{8} \\
   e_{8}
 \end{pmatrix} =   \begin{pmatrix}
  -2^{15}\cdot 4740272 \\
  -2^{18}\cdot 4756144\\
 -2^{21}\cdot 4787888 \\
   \end{pmatrix} 
    \end{equation*}
    \begin{equation} 
    \label{eq 2.21}
\displaystyle  \Leftrightarrow
    \begin{pmatrix}
   c_{8} \\
   d_{8} \\
   e_{8}
    \end{pmatrix} = \begin{pmatrix}
 -2^5\cdot 33626320\\
2^{10}\cdot 67570080\\
-2^{15}\cdot 38684032
  \end{pmatrix} 
  \end{equation}
Thus  :\\
\begin{align}
\label{eq 2.22}
\Gamma_{8}^{\left[2\atop {2\atop {2\atop{2\atop2 }}}\right]\times k}   = 
 496 \cdot 2^{4k}+496 \cdot 9525 \cdot 2^{3k} -2^5 \cdot 33626320 \cdot2^{2k}+2^{10}\cdot 67570080 \cdot 2^{k} -2^{15}\cdot 38684032 \\
 = 496\cdot [2^{4k}+9525\cdot2^{3k}-2169440\cdot2^{2k}+68115\cdot2^{11}\cdot2^{k}-9749\cdot2^{18}] \nonumber
\end{align}
To sum up we deduce from \eqref{eq 2.22},\eqref{eq 2.20}   :   \\

 \begin{equation}
 \label{eq 2.23}
  \begin{cases} 
 \displaystyle   \Gamma_{0}^{\left[2\atop {2\atop {2\atop{2\atop2 }}}\right]\times k}   = 1 \quad \text{for} \quad  k\geqslant 1 \\
\displaystyle    \Gamma_{1}^{\left[2\atop {2\atop {2\atop{2\atop2 }}}\right]\times k}   = 93 \quad \text{for} \quad  k\geqslant 2 \\
\displaystyle   \Gamma_{2}^{\left[2\atop {2\atop {2\atop{2\atop2 }}}\right]\times k}  = 31\cdot2^{k+1}+6386 \quad \text{for} \quad k\geqslant 3\\
\displaystyle    \Gamma_{3}^{\left[2\atop {2\atop {2\atop{2\atop2 }}}\right]\times k}  = 6510\cdot2^{k}+364560 \quad \text{for} \quad k\geqslant 4\\
\displaystyle    \Gamma_{4}^{\left[2\atop {2\atop {2\atop{2\atop2 }}}\right]\times k}  = 620\cdot2^{2k}+448260\cdot2^{k}+15748000 \quad \text{for} \quad k\geqslant 5\\
  \displaystyle    \Gamma_{5}^{\left[2\atop {2\atop {2\atop{2\atop2 }}}\right]\times k}  =  65100\cdot2^{2k}+22654800\cdot2^{k}+250817280   \quad \text{for} \quad k\geqslant 6\\
   \displaystyle    \Gamma_{6}^{\left[2\atop {2\atop {2\atop{2\atop2 }}}\right]\times k}  =  1240\cdot 2^{3k}+b_{6}\cdot2^{2k} +c_{6}\cdot 2^{k}  + d_{6}   \quad \text{for  }   \quad  k\geqslant 7 \\
 \displaystyle    \Gamma_{7}^{\left[2\atop {2\atop {2\atop{2\atop2 }}}\right]\times k} = 115320\cdot 2^{3k}+b_{7}\cdot2^{2k}+c_{7}\cdot2^{k} +d_{7}  \quad \text{for  } \quad  k\geqslant 8  \\
 \displaystyle    \Gamma_{8}^{\left[2\atop {2\atop {2\atop{2\atop2 }}}\right]\times k} = 496\cdot [2^{4k}+9525\cdot2^{3k}-2169440\cdot2^{2k}+68115\cdot2^{11}\cdot2^{k}-9749\cdot2^{18}]  \quad  \text{for   } \quad  k\geqslant 9  \\
 \displaystyle     \Gamma_{9}^{\left[2\atop {2\atop {2\atop{2\atop2 }}}\right]\times k} =31248\cdot[2^{4k}-480\cdot2^{3k} +71680\cdot2^{2k}-3932160\cdot2^{k}+2^{26}]  \quad  \text{for   }\quad  k\geqslant 10   \\   
 \displaystyle   \Gamma_{10}^{\left[2\atop {2\atop {2\atop{2\atop2 }}}\right]\times k}  = 2^5\cdot[2^{5k}-992\cdot2^{4k}+317440\cdot2^{3k}-40632320\cdot2^{2k}+2080374784\cdot2^{k}-2^{35}]  \; \text{for} \; k\geqslant 10 \\
\end{cases}
    \end{equation}
\textbf{Step 8}
\newline
Computation of $b_{6},b_{7} $ in \eqref{eq 2.23}.\\
From \eqref{eq 2.15} and \eqref{eq 2.23} we get :\\
 \begin{equation*}
   \displaystyle 
    \begin{pmatrix}
   1 & 1 \\ 
    2^{4} & 2^{3}    \\
   2^{8} &  2^6  
   \end{pmatrix}\displaystyle  \begin{pmatrix}
   b_{6} \\
   b_{7} 
   \end{pmatrix} =   \begin{pmatrix}
  136354120 \\
  1122567040\\
 9488281600 
   \end{pmatrix} 
    \end{equation*}
    \begin{equation} 
    \label{eq 2.24}
\displaystyle  \Leftrightarrow
    \begin{pmatrix}
   b_{6} \\
   b_{7} 
    \end{pmatrix} = \begin{pmatrix}
 1240 \cdot 3199\\
115320 \cdot 1148\\
 \end{pmatrix} 
  \end{equation}
To sum up we deduce from \eqref{eq 2.23},\eqref{eq 2.24}  : \\
 \begin{equation}
 \label{eq 2.25}
  \begin{cases} 
 \displaystyle   \Gamma_{0}^{\left[2\atop {2\atop {2\atop{2\atop2 }}}\right]\times k}   = 1 \quad \text{for} \quad  k\geqslant 1 \\
\displaystyle    \Gamma_{1}^{\left[2\atop {2\atop {2\atop{2\atop2 }}}\right]\times k}   = 93 \quad \text{for} \quad  k\geqslant 2 \\
\displaystyle   \Gamma_{2}^{\left[2\atop {2\atop {2\atop{2\atop2 }}}\right]\times k}  = 31\cdot2^{k+1}+6386 \quad \text{for} \quad k\geqslant 3\\
\displaystyle    \Gamma_{3}^{\left[2\atop {2\atop {2\atop{2\atop2 }}}\right]\times k}  = 6510\cdot2^{k}+364560 \quad \text{for} \quad k\geqslant 4\\
\displaystyle    \Gamma_{4}^{\left[2\atop {2\atop {2\atop{2\atop2 }}}\right]\times k}  = 620\cdot2^{2k}+448260\cdot2^{k}+15748000 \quad \text{for} \quad k\geqslant 5\\
  \displaystyle    \Gamma_{5}^{\left[2\atop {2\atop {2\atop{2\atop2 }}}\right]\times k}  =  65100\cdot2^{2k}+22654800\cdot2^{k}+250817280   \quad \text{for} \quad k\geqslant 6\\
   \displaystyle    \Gamma_{6}^{\left[2\atop {2\atop {2\atop{2\atop2 }}}\right]\times k}  =  1240\cdot 2^{3k}+ 1240 \cdot 3199 \cdot 2^{2k} +c_{6}\cdot 2^{k}  + d_{6}   \quad \text{for  }   \quad  k\geqslant 7 \\
 \displaystyle    \Gamma_{7}^{\left[2\atop {2\atop {2\atop{2\atop2 }}}\right]\times k} = 115320\cdot 2^{3k}+115320 \cdot 1148 \cdot 2^{2k}+c_{7}\cdot2^{k} +d_{7}  \quad \text{for  } \quad  k\geqslant 8  \\
 \displaystyle    \Gamma_{8}^{\left[2\atop {2\atop {2\atop{2\atop2 }}}\right]\times k} = 496\cdot [2^{4k}+9525\cdot2^{3k}-2169440\cdot2^{2k}+68115\cdot2^{11}\cdot2^{k}-9749\cdot2^{18}]  \quad  \text{for   } \quad  k\geqslant 9  \\
 \displaystyle     \Gamma_{9}^{\left[2\atop {2\atop {2\atop{2\atop2 }}}\right]\times k} =31248\cdot[2^{4k}-480\cdot2^{3k} +71680\cdot2^{2k}-3932160\cdot2^{k}+2^{26}]  \quad  \text{for   }\quad  k\geqslant 10   \\   
 \displaystyle   \Gamma_{10}^{\left[2\atop {2\atop {2\atop{2\atop2 }}}\right]\times k}  = 2^5\cdot[2^{5k}-992\cdot2^{4k}+317440\cdot2^{3k}-40632320\cdot2^{2k}+2080374784\cdot2^{k}-2^{35}]  \; \text{for} \; k\geqslant 10 \\
\end{cases}
    \end{equation}

\textbf{Step 9}
\newline
Computation of $c_{7},d_{7} $ in \eqref{eq 2.25}.\\
From \eqref{eq 2.13} we obtain :\\
 \begin{equation*}
   \displaystyle 
    \begin{pmatrix}
   32 & 1 \\ 
    64 & 1   
     \end{pmatrix}\displaystyle  \begin{pmatrix}
   c_{7} \\
   d_{7} 
   \end{pmatrix} =   \begin{pmatrix}
  -2^{15}\cdot 4252425 \\
  -2^{17}\cdot 4367745
 \end{pmatrix} 
    \end{equation*}
    \begin{equation} 
    \label{eq 2.26}
\displaystyle  \Leftrightarrow
    \begin{pmatrix}
   c_{7} \\
   d_{7} 
    \end{pmatrix} = \begin{pmatrix}
 -2^{10}\cdot 13218555 = 115320\cdot (-2^7\cdot 917)\\
2^{15}\cdot 8966130 = 115320\cdot (2^{12}\cdot 622 )
 \end{pmatrix} 
  \end{equation}
To sum up we deduce from \eqref{eq 2.25},\eqref{eq 2.26}  : \\
 \begin{equation}
 \label{eq 2.27}
  \begin{cases} 
 \displaystyle   \Gamma_{0}^{\left[2\atop {2\atop {2\atop{2\atop2 }}}\right]\times k}   = 1 \quad \text{for} \quad  k\geqslant 1 \\
\displaystyle    \Gamma_{1}^{\left[2\atop {2\atop {2\atop{2\atop2 }}}\right]\times k}   = 93 \quad \text{for} \quad  k\geqslant 2 \\
\displaystyle   \Gamma_{2}^{\left[2\atop {2\atop {2\atop{2\atop2 }}}\right]\times k}  = 31\cdot2^{k+1}+6386 \quad \text{for} \quad k\geqslant 3\\
\displaystyle    \Gamma_{3}^{\left[2\atop {2\atop {2\atop{2\atop2 }}}\right]\times k}  = 6510\cdot2^{k}+364560 \quad \text{for} \quad k\geqslant 4\\
\displaystyle    \Gamma_{4}^{\left[2\atop {2\atop {2\atop{2\atop2 }}}\right]\times k}  = 620\cdot2^{2k}+448260\cdot2^{k}+15748000 \quad \text{for} \quad k\geqslant 5\\
  \displaystyle    \Gamma_{5}^{\left[2\atop {2\atop {2\atop{2\atop2 }}}\right]\times k}  =  65100\cdot2^{2k}+22654800\cdot2^{k}+250817280   \quad \text{for} \quad k\geqslant 6\\
   \displaystyle    \Gamma_{6}^{\left[2\atop {2\atop {2\atop{2\atop2 }}}\right]\times k}  =  1240\cdot 2^{3k}+ 1240 \cdot 3199 \cdot 2^{2k} +c_{6}\cdot 2^{k}  + d_{6}   \quad \text{for  }   \quad  k\geqslant 7 \\
 \displaystyle    \Gamma_{7}^{\left[2\atop {2\atop {2\atop{2\atop2 }}}\right]\times k} = 115320\cdot [ 2^{3k}+1148 \cdot 2^{2k}-2^7 \cdot 917\cdot2^{k} + 2^{12}\cdot 622]  \quad \text{for  } \quad  k\geqslant 8  \\
 \displaystyle    \Gamma_{8}^{\left[2\atop {2\atop {2\atop{2\atop2 }}}\right]\times k} = 496\cdot [2^{4k}+9525\cdot2^{3k}-2169440\cdot2^{2k}+68115\cdot2^{11}\cdot2^{k}-9749\cdot2^{18}]  \quad  \text{for   } \quad  k\geqslant 9  \\
 \displaystyle     \Gamma_{9}^{\left[2\atop {2\atop {2\atop{2\atop2 }}}\right]\times k} =31248\cdot[2^{4k}-480\cdot2^{3k} +71680\cdot2^{2k}-3932160\cdot2^{k}+2^{26}]  \quad  \text{for   }\quad  k\geqslant 10   \\   
 \displaystyle   \Gamma_{10}^{\left[2\atop {2\atop {2\atop{2\atop2 }}}\right]\times k}  = 2^5\cdot[2^{5k}-992\cdot2^{4k}+317440\cdot2^{3k}-40632320\cdot2^{2k}+2080374784\cdot2^{k}-2^{35}]  \; \text{for} \; k\geqslant 10 \\
\end{cases}
    \end{equation}
\textbf{Step 10}
\newline
Computation of $c_{6} $ in \eqref{eq 2.27}.\\
From \eqref{eq 2.15} and \eqref{eq 2.27} we get :\\
\begin{align}
\label{eq 2.28}
2^{k}\cdot \big( 62+6510+448260+22654800+c_{6}+(-115320\cdot 2^7\cdot 917)+\\
496\cdot 68115 \cdot 2^{11}+(-31248 \cdot 3932160) +2^5 \cdot 2080374784 \big) = 0 \nonumber \\
\Leftrightarrow \quad c_{6} = 2^{10}\cdot 606515 = 1240Ê\cdot 3913\cdot 2^7 \nonumber
\end{align}

To sum up we deduce from \eqref{eq 2.28},\eqref{eq 2.27} :  \\
 \begin{equation}
 \label{eq 2.29}
  \begin{cases} 
 \displaystyle   \Gamma_{0}^{\left[2\atop {2\atop {2\atop{2\atop2 }}}\right]\times k}   = 1 \quad \text{for} \quad  k\geqslant 1 \\
\displaystyle    \Gamma_{1}^{\left[2\atop {2\atop {2\atop{2\atop2 }}}\right]\times k}   = 93 \quad \text{for} \quad  k\geqslant 2 \\
\displaystyle   \Gamma_{2}^{\left[2\atop {2\atop {2\atop{2\atop2 }}}\right]\times k}  = 31\cdot2^{k+1}+6386 \quad \text{for} \quad k\geqslant 3\\
\displaystyle    \Gamma_{3}^{\left[2\atop {2\atop {2\atop{2\atop2 }}}\right]\times k}  = 6510\cdot2^{k}+364560 \quad \text{for} \quad k\geqslant 4\\
\displaystyle    \Gamma_{4}^{\left[2\atop {2\atop {2\atop{2\atop2 }}}\right]\times k}  = 620\cdot2^{2k}+448260\cdot2^{k}+15748000 \quad \text{for} \quad k\geqslant 5\\
  \displaystyle    \Gamma_{5}^{\left[2\atop {2\atop {2\atop{2\atop2 }}}\right]\times k}  =  65100\cdot2^{2k}+22654800\cdot2^{k}+250817280   \quad \text{for} \quad k\geqslant 6\\
   \displaystyle    \Gamma_{6}^{\left[2\atop {2\atop {2\atop{2\atop2 }}}\right]\times k}  =  1240\cdot 2^{3k}+ 1240 \cdot 3199 \cdot 2^{2k} 
   +1240Ê\cdot 3913\cdot 2^7 \cdot 2^{k}  + d_{6}   \quad \text{for  }   \quad  k\geqslant 7 \\
 \displaystyle    \Gamma_{7}^{\left[2\atop {2\atop {2\atop{2\atop2 }}}\right]\times k} = 115320\cdot [ 2^{3k}+1148 \cdot 2^{2k}-2^7 \cdot 917\cdot2^{k} + 2^{12}\cdot 622]  \quad \text{for  } \quad  k\geqslant 8  \\
 \displaystyle    \Gamma_{8}^{\left[2\atop {2\atop {2\atop{2\atop2 }}}\right]\times k} = 496\cdot [2^{4k}+9525\cdot2^{3k}-2169440\cdot2^{2k}+68115\cdot2^{11}\cdot2^{k}-9749\cdot2^{18}]  \quad  \text{for   } \quad  k\geqslant 9  \\
 \displaystyle     \Gamma_{9}^{\left[2\atop {2\atop {2\atop{2\atop2 }}}\right]\times k} =31248\cdot[2^{4k}-480\cdot2^{3k} +71680\cdot2^{2k}-3932160\cdot2^{k}+2^{26}]  \quad  \text{for   }\quad  k\geqslant 10   \\   
 \displaystyle   \Gamma_{10}^{\left[2\atop {2\atop {2\atop{2\atop2 }}}\right]\times k}  = 2^5\cdot[2^{5k}-992\cdot2^{4k}+317440\cdot2^{3k}-40632320\cdot2^{2k}+2080374784\cdot2^{k}-2^{35}]  \; \text{for} \; k\geqslant 10 \\
\end{cases}
    \end{equation}
\textbf{Step 11}
\newline
Computation of $d_{6} $ in \eqref{eq 2.29}  :  \\
From \eqref{eq 2.13} we deduce : \\
\begin{align}
\label{eq 2.30}
 \Gamma_{6}^{\left[2\atop {2\atop {2\atop{2\atop2 }}}\right]\times 5}  = 0 \\
 \Leftrightarrow \quad 
  1240\cdot 2^{15}+ 1240 \cdot 3199 \cdot 2^{10} +1240Ê\cdot 3913\cdot 2^7 \cdot 2^{5}  + d_{6} = 0 \nonumber \\
  \Leftrightarrow \quad    d_{6} = - 1240\cdot [ 2^{15}+ 3199 \cdot 2^{10} + 3913 \cdot 2^{12}] = -1240 \cdot 18883\cdot 2^{10} \nonumber
\end{align}
and Theorem \ref{thm 2.1}  is proved.

\end{proof}

\end{document}